\documentclass[reqno,10pt,centertags,draft]{amsart}
\usepackage{amsmath,amsthm,amscd,amssymb,latexsym,upref}
\date{\today}
\newcommand{\bbD}{{\mathbb{D}}}

\newcommand{\bbZ}{{\mathbb{Z}}}

\newcommand{\bbT}{{\mathbb{T}}}

\newcommand{\cH}{{\mathcal{H}}}

\allowdisplaybreaks
\numberwithin{equation}{section}
\newtheorem{theorem}{Theorem}[section]

\newtheorem{lemma}[theorem]{Lemma}

\newtheorem{corollary}[theorem]{Corollary}

\theoremstyle{definition}
\newtheorem{definition}[theorem]{Definition}
\newtheorem{remark}[theorem]{Remark}

\begin{document}

\title[]{On a new asymptotic problem in
the scattering setting}
\author[]{ F. Peherstorfer, A. Volberg,
and P. Yuditskii}

\thanks{Partially supported by NSF grant DMS-0200713 and
the Austrian Science Found FWF, project number: P16390--N04 }

\address{Institute for Analysis, Johannes Kepler University Linz,
A-4040 Linz, Austria}

\email{Franz.Peherstorfer@jku.at}

\address{
Department of Mathematics, Michigan State University, East
Lansing, Michigan 48824, USA}
\email{volberg@math.msu.edu}

\address{
Department of Mathematics, Bar Ilan University, Tel Aviv, Israel}
\email{yuditski@macs.biu.ac.il}

\date{\today}

\begin{abstract} In recent works we considered an asymptotic problem for orthogonal polynomials
when a Szeg\"o measure on the unit circumference is perturbed by 
an arbitrary Blaschke sequence of point masses outside the unit
disk. In the current work we consider a similar problem in the
scattering setting.
\end{abstract}

\maketitle

The goal of this work is to consider a new asymptotic problem; the
related problems in the spectral setting were solved in
\cite{NVYu}, \cite{PVYu}.

  With a given Szeg\"o contractive function $R$ on the unite
  circle $\bbT$
  \begin{equation}\label{is1}
  |R(t)|\le 1,\quad \log(1-|R|)\in L^1
  \end{equation}
  and a positive measure $\nu$ supported on the  Blaschke set
  \begin{equation}\label{is2}
  Z=\{\zeta_k:\sum(1-|\zeta_k|)<\infty\},\quad
  \nu(\zeta_k)=\nu_{k}>0,
  \end{equation}
  we associate the scalar product
  \begin{equation}\label{is3}
\langle D f,f\rangle=\langle (I-\Gamma^*\Gamma)f,f\rangle +
\sum|f(\zeta_k)|^2\nu_{k}.
\end{equation}
Here $\langle\cdot,\cdot\rangle$ is the inner standard product in
$L^2$ on $\bbT$ with respect to the Lebesgue measure, $\Gamma$ is
the Hankel operator with the symbol $R$,
$$
\Gamma f=\Gamma_{R} f=P_- (R f),
$$
acting from the Hardy space $H^2$ into its orthogonal complement
$H^2_-=L^2\ominus H^2$, $P_-$ is the ortho--projection onto
$H^2_-$. Let us point out that we even do not require that the
measure $\nu$ is finite, thus the scalar product, correspondingly
the unbounded operator $D$, are defined initially on the
$H^2$--functions that equal zero at all points of $Z$ except for a
finite number of them. (In this place we use the Blaschke
condition \eqref{is2}).

Our generalization deals with the presence of the measure $\nu$.
Without it the scalar product plays the key role in the classical
now solution of the Nehari problem by Adamyan, Arov and Krein,
\cite{AAK1}, \cite{AAK2}. For the new point of view on this
subject see \cite{TTh}. The Nehari problem is also known as the
generalized Schur problem. Concerning relation of the Schur
problem with the Theory of Orthogonal Polynomials on the Unit
Circle, CMV matrices and so on, see \cite{S1}, \cite{S2}.

For shortness we denote the collection of data by
$$
\alpha:=\{R,\nu\}
$$
and then use the notation  $H^2(\alpha)$ for the closure of
admissible functions $f$ from $H^2$ with respect to the metric
 \eqref{is3} with $D=D(\alpha)$. The condition \eqref{is1}
guaranties that the point evaluation functional (for all $\zeta_0$
in the unit disk $\bbD$)
$$
f\mapsto f(\zeta_0)
$$
 is bounded in $H^2(\alpha)$. Let $k^{\alpha}$ be the
 reproducing kernel of this space with respect to the origin and let
 $$
K^{\alpha}=\frac{k^{\alpha}}{\Vert k^{\alpha}\Vert}.
$$
It is almost evident that the system
\begin{equation}\label{is4}
\{e_n(\zeta;\alpha)\},\quad e_n(\zeta;\alpha):=\zeta^n
K^{\alpha_{n}}(\zeta),
\end{equation}
where
$$
\alpha_{n}=\{\zeta^n R(\zeta),|\zeta|^{2n}\nu(\zeta)\}
$$
forms an orthonormal basis in $H^2(\alpha)$. We claim that
asymptotically this system behaves as the standard basis system in
$H^2$, in particulary, that
\begin{equation}\label{is5}
K^{\alpha_{n}}(0)\to 1, \quad n\to\infty.
\end{equation}

 We follow the line of proof that was suggested in \cite{PYu} and
then improved in \cite{VYu} and \cite{KPVYu}. Actually, the
general idea is very simple. There are two natural steps in
approximation of the given data by ``regular" ones. First, to
substitute the given measure $\nu$ by a finitely supported
$\nu^N$. Second, to substitute $R$ by $\rho R$ with $0<\rho<1$.
Then the corresponding data $\alpha^{N,\rho}$ produce the metric
$D(\alpha^{N,\rho})$ which is equivalent to the standard metric in
$H^2$ and it is a fairly easy task to prove \eqref{is5} for such
data. Further, due to $D(\alpha^{N})\le D(\alpha)\le
D(\alpha^{\rho})$ we have the evident estimations
$$
K^{\alpha^{N}}(0)\ge K^{\alpha}(0)\ge K^{\alpha^{\rho}}(0).
$$
And the key point is a certain duality principle, see Corollary
\ref{cor1.6}, that will allow us to use the left or right side
estimation whenever it is convenient for us.

\section{The duality}
\subsection{The space $L^2(\alpha)$}
We define the outer function $T_e$ by
$$
|T_e|^2=1-|R|^2,\ \  T_e(0)>0.
$$
Consider  the scalar product
\begin{equation}
\left\langle \begin{bmatrix} 1&\bar R\\
                             R& 1
                             \end{bmatrix}
                             \begin{bmatrix} f_1\\
                             f_2
                             \end{bmatrix},
                             \begin{bmatrix} f_1\\
                             f_2
                             \end{bmatrix}
 \right\rangle.
\end{equation}
By $L^2_{R}$ we denote the closure of  vectors of the form
\begin{equation}\label{db23}
\begin{bmatrix} f\\
                          -P_- R f
                             \end{bmatrix}, \quad f\in L^2.
\end{equation}
with respect to the above metric.

\begin{lemma} We have
 \begin{equation}\label{db29}
  L^2_{R}=
\left\{
\begin{bmatrix} f_1\\f_2
\end{bmatrix}: T_ef_1\in L^2, \ \bar T_e f_2 \in H^2_-,
  \
 R f_1+ f_2\in H^2
                             \right\}.
\end{equation}
Moreover, the first component $f_1$ determines the second
component $f_2$ uniquely.
\end{lemma}

\begin{proof}
As usually we can find the closure as the second orthogonal
complement. This proves all listed properties of the vector
$\begin{bmatrix} f_1\\f_2
\end{bmatrix}$.

To prove the second claim let us mention that the spaces
$$
\frac{1}{T_e}H^2\dotplus H^2_- \quad \text{and} \quad
\frac{1}{\bar T_e}H^2_-\dotplus H^2
$$
form direct sums (due to the maximum principle in the Smirnov
class: $\frac{g_1}{g_2}\in L^2$ with $g_{1,2}\in H^\infty$, $g_2$
is outer, implies $\frac{g_1}{g_2}\in H^2$). Thus the
decomposition
$$
R f_1=-\frac{1}{\bar T_e}(\bar T_e f_2)+h, \quad h\in H^2,
$$
 is unique.
\end{proof}

\begin{definition}
The space $L^2(\alpha)$ is formed by functions $f$ that are
defined on $\bbT\cup Z$ and such that $f\vert_{\bbT}=f_1$,
$\begin{bmatrix} f_1\\f_2
\end{bmatrix}\in L^2_R$, and $f(\zeta_k)=f_k$,
$\{f_k\}\in L^2_{\nu}$,  with the norm
$$
\Vert f\Vert^2_{L^2(\alpha)}=\left\Vert\begin{bmatrix} f_1\\f_2
\end{bmatrix}\right\Vert^2_{L^2_R}+
\Vert \{f_k\}\Vert^2_{L^2_{\nu}}.
$$
In other words $L^2(\alpha)=L^2_R\oplus L^2_\nu$.
\end{definition}

\subsection{The Hardy spaces $\check{H^2}(\alpha)$ and
$\hat{H^2}(\alpha)$} The first space $\check H^2(\alpha)$ is
defined as the closure of functions $f\in H^2$, that equal zero at
all points of $Z$ except for a finite number of them, in
$L^2(\alpha)$ (a function from $H^2$ is naturally defined on
$\bbT\cup Z$).
\begin{lemma}\label{isl1.3}
If $f\in \check H^2(\alpha)$ then
\begin{equation}\label{is1.4}
g:= T_e f\vert_{\bbT}\in H^2\quad \text{and} \quad
    f(\zeta_k)=\frac{g(\zeta_k)}{T_e(\zeta_k)}.
\end{equation}
\end{lemma}

\begin{definition}\label{def1.4}
A function $f\in L^2(\alpha)$ is in $\hat H^2(\alpha)$ if
conditions \eqref{is1.4} hold.
\end{definition}
Thus we have evidently
$\check{H^2}(\alpha)\subset\hat{H^2}(\alpha)$ but they do not
necessarily coincide. Of course, in a regular case they are the
same, say, when topologically both spaces are equivalent to the
standard $H^2$.

The both spaces have the reproducing kernel basis  systems:
\begin{equation}\label{is1.5}
\{\zeta^n \check K^{\alpha_{n}}(\zeta)\},\quad \{\zeta^n \hat
K^{\alpha_{n}}(\zeta)\},
\end{equation}
where $n\in \bbZ_+$. Let us point out that both systems being
extended to all $n\in \bbZ$ are basises in $L^2(\alpha)$.

\subsection{The dual $L^2$--space} Our goal is to describe
the orhogonal complement of, say, $\check H^2(\alpha)$ in
$L^2(\alpha)$. We will see that under a certain unitary map from
the given $L^2(\alpha)$ to a similar $L^2$--space such orthogonal
complement transforms into $\hat H^2$--subspace of the target
space.

Now we define formally this dual space.
 We define the Blaschke product
$$
B(\zeta)=\prod_{\zeta_k\in
Z}\frac{\zeta_k-\zeta}{1-\bar\zeta_k\zeta}\frac{|\zeta_k|}{\zeta_k}
$$
and the function $T=\frac{T_e}{B}$.

Put
\begin{equation}\label{is1.6}
    R^\tau(\bar t)=-\frac{R(t)}{T_e(t)}\overline{ T(t)}, \ t\in\bbT.
\end{equation}
Note that
\begin{equation}\label{is1.7}
    \begin{bmatrix} 1&\overline{R^\tau}\\
    R^\tau& 1
\end{bmatrix}(\bar t)=
\begin{bmatrix} T&0\\
    0& \bar T_e
\end{bmatrix}( t)
\begin{bmatrix} 1&\overline{R}\\
    R & 1
\end{bmatrix}^{-1}(t)
\begin{bmatrix} \bar T&0\\
    0&  T_e
\end{bmatrix}( t),
\end{equation}
and, also, due to $|R^\tau(\bar t)|=|R( t)|$, we have $T_e^\tau(
t)=\overline{T_e(\bar t)}$.

We define the measure $\nu^\tau$ supported on $Z^\tau:=\bar Z$ by
\begin{equation}\label{is1.8}
    \nu^\tau(\bar\zeta_k)
    \nu(\zeta_k)=\left|\left(\frac{1}{T}\right)'(\zeta_k)\right|^2.
\end{equation}
We have $B^\tau(t)=\overline{B(\bar t)}$, and thus
$T^\tau(t)=\overline{T(\bar t)}$.

Finally, the map $\tau: L^2(\alpha)\to L^2(\alpha^\tau)$ is
defined by
\begin{equation}\label{is1.9}
\begin{bmatrix} 1&\overline{R}\\
    R & 1
\end{bmatrix}(t)
\begin{bmatrix} f_1\\
    f_2
\end{bmatrix}(t)=
\begin{bmatrix} \bar T&0\\
    0&  T_e
\end{bmatrix}( t)\bar t
\begin{bmatrix} f_1^\tau\\
    f_2^\tau
\end{bmatrix}( \bar t)
\end{equation}
on $L^2_R$ component, and by
\begin{equation}\label{is1.10}
    f^\tau(\bar\zeta_k)=-\overline{
    \left(\frac{1}{T}
    \right)'(\zeta_k)}f(\zeta_k)\nu_k,
\end{equation}
on $L^2_\nu$, so that
\begin{equation}\label{is1.11}
\tau f\vert_{\bbT}=f_1^\tau, \quad (\tau f)(\bar\zeta_k)=
f^\tau(\bar\zeta_k).
\end{equation}

It is an easy task to check correctness of the definition, as well
as the fact that the map is an involution.

\subsection{The duality Theorem}

\begin{theorem}$\tau$ maps unitary $L^2(\alpha)\ominus \check
H^2(\alpha)$ onto $\hat H^2(\alpha^\tau)$.
\end{theorem}

\begin{proof}
Let $f$, with the components
\begin{equation*}
\begin{bmatrix}
f_1\\f_2
\end{bmatrix}
\in L^2_{R} ,\quad \{f_k\}\in L^2_{\nu},
\end{equation*}
be a vector from the orthogonal complement to $\check
H^2(\alpha)$. By orthogonality to the vectors of the form $Bh$,
$h\in H^2$, we get
\begin{equation}\label{is1.12}
\langle f, Bh\rangle_{L^2(\alpha)} =
 \left\langle
\begin{bmatrix}
\bar T \bar t f_1^\tau(\bar t)\\ T_e\bar t f_2^\tau(\bar t)
\end{bmatrix},
\begin{bmatrix}
Bh\\ -P_- R Bh
\end{bmatrix}
\right\rangle= \langle \bar T_e\bar t f_1^\tau(\bar t),  h\rangle
 =0.
\end{equation}
That is $T_e^\tau f_1^\tau\in H^2$.

Substituting in the scalar product the function
$\frac{B(t)}{t-\zeta_k}$ we get
\begin{equation*}
\left\langle \bar T_e\bar t f_1^\tau(\bar t),
\frac{1}{t-\zeta_k}\right\rangle
+f_k\overline{B'(\zeta_k)}\nu_k=0,
\end{equation*}
or,
\begin{equation*}
(T_e^\tau f_1^\tau)(\bar\zeta_k)+f_k\overline{B'(\zeta_k)}\nu_k=0.
\end{equation*}
By Definition \ref{def1.4} the target vector $\tau f$ is in $\hat
  H^2(\alpha^\tau)$.

  Conversely, an arbitrary vector of this form is orthogonal to
  $Bh$, $h\in H^2$ and $\frac{B(t)}{t-\zeta_k}$, $\forall k$, and
  such vectors are complete in $\check H^2(\alpha)$.

\end{proof}

Note that all definitions \eqref{is1.6} ... \eqref{is1.11} were
given to suite the proof of this theorem.

\begin{corollary}\label{cor1.6} In the above notations
\begin{equation}\label{is1.13}
    T(0)\check K^{\alpha_{-1}}(0)\hat K^{\alpha^{\tau}}(0)=
    1.
\end{equation}
\end{corollary}
\begin{proof}
By the theorem we conclude that
\begin{equation}\label{is1.14}
\tau (\zeta^{-1} \check K^{\alpha_{-1}}(\zeta))= \hat
K^{\alpha^{\tau}}(\zeta).
\end{equation}
 Since $\check K^{\alpha_{-1}}(\zeta)\check K^{\alpha_{-1}}(0)$
is the reproducing kernel of $\check H^2(\alpha_{-1})$, we have
\begin{equation}\label{is1.15}
\langle \zeta^{-1}\check K^{\alpha_{-1}}(\zeta), \zeta^{-1}
B(\zeta)\rangle _{L^2(\alpha)}= \langle \check
K^{\alpha_{-1}}(\zeta),  B(\zeta)\rangle _{L^2(\alpha_{-1})}=
\frac{B(0)}{K^{\alpha_{-1}}(0)}.
\end{equation}
On the other hand, using \eqref{is1.14}, in the same way as in
\eqref{is1.12}, we get for this scalar product
\begin{equation*}
    \langle \bar T\bar t K^{\alpha^\tau}(\bar t), \bar t B\rangle
    =T_e(0)K^{\alpha^\tau}(0).
\end{equation*}
In combination with \eqref{is1.15} it gives us \eqref{is1.13}.
\end{proof}

\section{Asymptotics}
\begin{proof}[Proof of \eqref{is5}] We have
\begin{equation}\label{is2.1}
\begin{split}
    \check K^{\alpha_n}(0) &\le  \check K^{\alpha^N_n}(0)=\frac{B^{(N)}(0)}{T_e(0)}
    \frac{1}{\hat K^{\alpha^N_{-n-1}}(0)}\\
    &\le
\frac{B^{(N)}(0)}{T_e(0)}
    \frac{1}{\hat K^{\alpha^{N,\rho}_{-n-1}}(0)}
    =
    \frac{T_e^{\rho}(0)}{T_e(0)}
    \check K^{\alpha^{N,\rho}_{n}}(0).
    \end{split}
\end{equation}
And from the other side
\begin{equation}\label{is2.2}
\begin{split}
    \check K^{\alpha_n}(0) &\ge  \check K^{\alpha^\rho_n}(0)=\frac{B(0)}{T^\rho_e(0)}
    \frac{1}{\hat K^{\alpha^\rho_{-n-1}}(0)}\\
    &\ge
\frac{B(0)}{T^\rho_e(0)}
    \frac{1}{\hat K^{\alpha^{N,\rho}_{-n-1}}(0)}
    =\frac{B(0)}{B^{(N)}(0)}
    \check K^{\alpha^{N,\rho}_{n}}(0).
    \end{split}
\end{equation}

Passing to the limit in \eqref{is2.1} and \eqref{is2.2} and using
arbitrerness of $\rho$ and $N$ we get \eqref{is5}.
\end{proof}

\begin{remark}
Assume that a space $\cH$, $\check
H^2(\alpha)\subset\cH\subset\hat H^2(\alpha)$ has the following
shift invariant property: $\check
H^2(\alpha_n)\subset\cH^{(n)}\subset\hat H^2(\alpha_n)$ for the
defined inductively
\begin{equation}\label{r}
\cH^{(n)}_0=\zeta \cH^{(n+1)},\ \ \cH^{(n)}_0:=\{f\in \cH^{(n)}:
f(0)=0\},
\end{equation}
with initial $\cH^{(0)}=\cH$. Then, naturally, there exists the
reproducing kernel based orthonormal basis in $\cH$ which also has
property \eqref{is5}.

\end{remark}

\bibliographystyle{amsplain}

\end{document}